\documentclass{svproc}
\usepackage{url}
\usepackage{graphicx}
\usepackage{mathrsfs}
\usepackage{multirow}
\usepackage{bm}
\usepackage{amsmath}
\usepackage{amsfonts}
\usepackage{amssymb}

\usepackage{subcaption}

\newcommand{\R}{\mathbb{R}}

\begin{document}
\mainmatter              
\title{Learning the solution operator of a nonlinear parabolic equation using physics informed deep operator network}
\titlerunning{Learning the solution operator using PI-DeepONet}  
%
\author{
Daniel \v{S}ev\v{c}ovi\v{c}\inst{1}
\and 
Cyril Izuchukwu Udeani\inst{1}
}
\authorrunning{Daniel \v{S}ev\v{c}ovi\v{c} and Cyril Izuchukwu Udeani} 
%
\tocauthor{Daniel \v{S}ev\v{c}ovi\v{c}, Cyril Izuchukwu Udeani}
\institute{
Comenius University in Bratislava, Mlynsk\'a dolina, 84248 Bratislava, Slovakia, \\
\email{sevcovic@fmph.uniba.sk},\\ 
\texttt{www.iam.fmph.uniba.sk/institute/sevcovic}
}

\maketitle              

\begin{abstract}
This study focuses on addressing the challenges of solving analytically intractable differential equations that arise in scientific and engineering fields such as Hamilton-Jacobi-Bellman. Traditional numerical methods and neural network approaches for solving such equations often require independent simulation or retraining when the underlying parameters change. To overcome this, this study employs a physics-informed DeepONet (PI-DeepONet) to approximate the solution operator of a nonlinear parabolic equation. PI-DeepONet integrates known physics into a deep neural network, which learns the solution of the PDE.

\keywords{Deep learning,  PI-DeepONet, Nonlinear parabolic equation}
\end{abstract}
\section{Introduction}
It is well-known that various nonlinear parabolic equations arise from various applied problems in industry. However, most of these differential equations are analytically intractable. Classical methods, such as the finite volume method, the finite difference method, and spectral methods, have been widely used to solve such equations. The corresponding finite-dimensional algebraic systems are often solved by iterative methods. Although these methods are efficient and well-studied, they require a lot of memory space and time, leading to high computational costs. Furthermore, a slight change in the input parameter leads to a new numerical simulation. To overcome these challenges, many researchers have replaced traditional numerical discretization methods with artificial neural networks (ANNs) to approximate the PDE solution. Recently, deep neural networks (DNNs) have been widely used to solve classical applied mathematical problems, including PDEs, utilizing machine learning and artificial intelligence approaches \cite{sirignano2018dgm}. Due to significant nonlinearities, convection dominance, or shocks, some PDEs are difficult to solve using standard numerical approaches. To this end, deep learning has recently emerged as a new paradigm of scientific computing thanks to the universal approximation theorem and the great expressivity of neural networks \cite{chen1995universal}.
Recent studies have shown that deep learning is a promising method for building metamodels for fast predictions of dynamic systems. In particular, neural networks (NNs) have been shown to represent the underlying nonlinear input-output relationship in complex systems. In an attempt to approximate the solution of PDEs, one can employ the deep Galerkin method \cite{sirignano2018dgm} involving DNNs to solve nonlinear PDEs. More recently, Lu et al. \cite{raissi2019physics} introduced an efficient technique called physics-informed neural networks (PINN) to approximate the solution of PDEs.   Although PINNs are faster than traditional numerical methods, they also have some limitations; e.g., a slight change in the underlying parameters could result in the retraining of the model. 
To overcome the shortcoming of PINNs, Lu et al. \cite{lu2021learning} further introduced the concept of DeepONet, which is an NN-based model that can learn linear and nonlinear PDE solution operators with a small generalization error via the universal approximation theorem for operators. DeepONet consists of two parts: a deep neural network that learns the solution of the PDE and an operator network that enforces the PDE at each iteration. The operator network acts as a constraint to ensure that the neural network outputs satisfy the underlying PDE. DeepONet maps input functions with infinite dimensions to output functions belonging to infinite-dimensional space. It can efficiently and accurately solve PDE with any initial and boundary conditions without retraining the network. PI-DeepONet approximates the PDE solution operator using two networks: one network that encodes the discrete input function space (branch net) and one that encodes the domain of the output functions (trunk net) (cf. \cite{lu2021learning}). It can effectively approximate the solution of different PDEs without requiring a large amount of training data by introducing a regularization mechanism that biases the output of DeepONet models to ensure physical consistency. PI-DeepONet can efficiently solve parametric linear and nonlinear PDEs compared to other variants of PINN since it can take source term parameters (including other parameters) as input variables. It can also break the curse of dimensionality in the input space, making it more suitable than other traditional approaches. 
Inspired by the above development and studies, we apply the PI-DeepONet approach for solving the following parabolic equation. 
\begin{equation}
\partial_\tau \varphi - \partial^2_x \alpha(\varphi) = g(\tau,x),\; (\tau, x)\in\Omega\equiv(0,T)\times (-L,L).
\label{eq_PDE-consdelta}
\end{equation}
For simplicity, we consider zero initial and boundary conditions for the solution $\varphi(\tau, x)$. Here $g$ is the source term. This model equation arises from the Hamilton-Jacobi-Bellman (HJB) equation describing the stochastic optimization problem (see \v{S}ev\v{c}ovi\v{c} and Kilianov\'a \cite{KilianovaSevcovicANZIAM} and \v{S}ev\v{c}ovi\v{c} and Udeani \cite{udeani2021application}).
The diffusion function $\alpha$ is the value function arising from a convex parametric optimization problem (see \v{S}ev\v{c}ovi\v{c} and Kilianov\'a \cite{KilianovaSevcovicANZIAM} and Kilianov\'a and Trnovsk\'a \cite{KilianovaTrnovska} for details).

\section{Methodology of PI-DeepONet}
\begin{figure}
\centering
\begin{subfigure}{.35\textwidth}
  \centering
     \includegraphics[width=\linewidth]{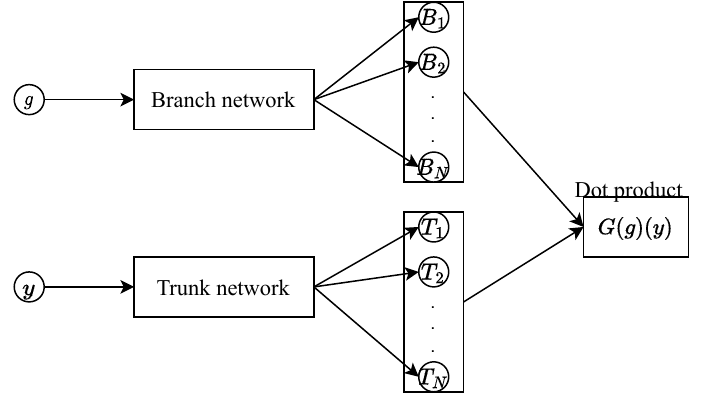}
  \caption{DeepONet}
  \label{fig_DeepONet}
\end{subfigure}%
\begin{subfigure}{.65\textwidth}
  \centering
    \includegraphics[width=\linewidth]{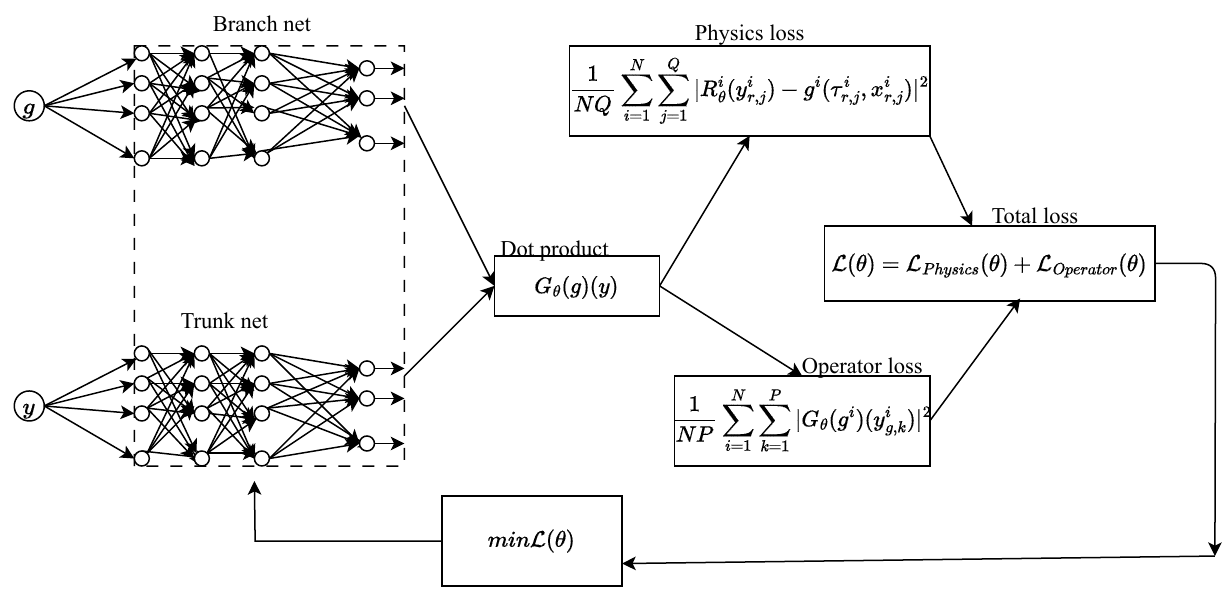}
  \caption{Physics informed DeepONet}
  \label{fig_piDeepONet}
\end{subfigure}
\caption{\small Schematics of DeepONet a), and physics informed DeepONet b)}
\label{fig_network}
\end{figure}

In this section, we introduce and discuss the methodology of PI-DeepONet. Consider the following equation:
\begin{equation}
    \mathcal{F} (g,\varphi) = 0,
\label{PDEop}
\end{equation}
where $\mathcal{F}$ is a differential operator for the governing PDE of some underlying physics laws, $g$ denotes its source term, and $\varphi$ is its solution. The differential equation (\ref{PDEop}) is assumed to have zero initial and boundary conditions. Note that the same idea can be applied to any initial and boundary conditions. Let $G: g \to G(g)$ be an operator between two infinite-dimensional function spaces where $g$ and $G(g)$ are two functions. This mapping is called the solution operator of equation (\ref{PDEop}), which can be evaluated at a random location $y$. In learning an operator in a more general setting, the inputs usually consist of two independent parts: the input function $g$ and the location variable (s) $y$. This learning can be done directly using traditional neural networks such as feedforward neural networks (FNN), recurrent neural networks (RNNs), convolutional neural networks (CNNs), or combining the two inputs as a single network input $(i.e., \{g, y\})$. Meanwhile, it is not necessarily advisable to directly use RNNs or CNNs since the input does not have a definite structure. Therefore, it is recommended to use FNNs as the baseline model. Furthermore, the DeepONet consists of branch and trunk nets. The branch net takes $g$ as the input function evaluated at a collection of fixed sensors $\{x_i\}^m_{i=1}$ and outputs a feature embedding of $q$ dimensions. The trunk net takes $y$ as input and also outputs a feature embedding of $q$ dimensions. Note that the dimensions of $y$ and $g$ need not be the same, indicating that $g$ and $y$ need not be treated as a single input like traditional NN. In general, the DeepONet network for learning an operator takes $g$ and $y$ as inputs and outputs $G(g)(y)$, which is obtained by taking the dot product of two subnetworks. The dot product of the outputs of the two subnets plays a crucial role in determining how well the learned solution operator aligns with the actual solution of the PDE. It measures the similarity or alignment between the two networks' outputs. This helps to improve the accuracy of the learned solution operator. Consequently, the PI-DeepONet is trained by minimizing the loss function $\mathcal{L}(\theta)$ (see \ref{totalloss}) over all the input-output triplets $\{g, y, G(g)(y)\}$,
where $\theta$ is the set of the weight matrix and the bias vector in the networks. The first goal is to find such an approximator $G_{\theta}(g)$, but thanks to the universal approximation theorem for operator \cite[Theorem 5]{chen1995universal}, which guarantees the existence of such function, i.e., $G_{\theta}(g)(y) \approx G(g)(y) = \varphi(y)\in \R$. The final objective is to find the best parameters that minimize the loss function (\ref{totalloss}) using suitable optimization techniques. The universal approximation theorem shows the stacked and unstacked DeepONet. The stacked network has one trunk net and $P$ stacked branch nets, whereas the unstacked network has one trunk net and one branch net, which are fully independently connected. For more details, see T. Chen and H. Chen \cite{chen1995universal}. Fig. \ref{fig_DeepONet} shows the schematics of an unstacked DeepONet. In this study, we use an unstacked DeepONet to solve a parametric parabolic equation arising from portfolio selection problems.

\section{Problem formulation}

To employ PI-DeepONet to solve the nonlinear parabolic equation (\ref{eq_PDE-consdelta}), we first define an operator that maps the input function to the PDE solution as
$G(g) =\varphi.$ The novelty of DeepONet is that it takes any arbitrary source term function as the input variables, making it more suitable than the PINN approach. Since $\varphi$ is also a function, we can evaluate it at some point, say $y$, to obtain $G(g)(y) =\varphi(y)$. In our application, $y =(\tau, x)$ denotes the point in the computational domain $\Omega$ where the network predicts the solution of the PDE (\ref{eq_PDE-consdelta}). In general, the branch (with $g$ as input function) and trunk (with $y$ as input variable) networks are given by $\mathcal{B}(\textbf{g}(\bf{\tilde{x}})) = \textbf{c}\cdot\sigma(\textbf{W}_{\mathcal{B}} \cdot g(\bf{\tilde{x}}) + \bf{b}_{\mathcal{B}})$ and $\mathcal{T}(y)=\sigma(\bf{W}_{\mathcal{T}}\cdot\bf{y} + \bf{b}_{\mathcal{T}}),$
respectively. Here, $\bf{\tilde{x}} = (\tau, x)$; $c$ is some positive constant; $\sigma$ is the activation function; $\bf{W}_{\mathcal{B}}$ and $\bf{W}_{\mathcal{T}}$ represent the weight matrices of branch and trunk networks, respectively; $\bf{b}_{\mathcal{B}}$ and $\bf{b}_{\mathcal{T}}$ represent the bias vector of branch and trunk networks, respectively. 

Now, letting $g^{i}, i=1,\dots, N$, be any given input function representing the source term in (\ref{eq_PDE-consdelta}), then equation (\ref{eq_PDE-consdelta}) becomes $g^i = \partial_{\tau} \varphi^i -\partial^2_x \alpha(\varphi^i)$.
According to \cite[Theorem 5]{chen1995universal}, there exists $G_{\theta}(g^i)$ such that $G_{\theta}(g^i)(y) \approx G(g^i)(y) = \varphi^i (y)$. For a fixed $i$, the approximator in the DeepONet solution operator is the dot product of the outputs of the branch and trunk networks, i.e.,
$G_{\theta}(g^i)(y) = \mathcal{B}(\textbf{g}(y)) \cdot \mathcal{T}(y).$ Hence, $g^i \approx \partial_{\tau} G_{\theta}(g^i)(y) -\partial^2_x \alpha(G_{\theta}(g^i)(y)).$
Therefore, the physics loss evaluated at the $Q$ collocation points in the interior of the domain is

\[
\mathcal{L}_{Physics} (\theta)= \frac{1}{NQ}\sum_{i=1}^N\sum_{j=1}^Q|R_{\theta}^i(y_{r, j}^i) - g^i ( x_{r, j}^i)|^2.
\]
Here, $R_{\theta}^i(y_{r, j}^i)= \partial_{\tau} G_{\theta}(g^i)(y_{r, j}^i) -\partial^2_x \alpha(G_{\theta}(g^i)(y_{r, j}^i))$ represents the residual that satisfies the underlying PDE, and $y_{r, j}^i = (\tau_{r, j}^i, x_{r,j}^i)$ denotes the collocation points where the PDE is evaluated. Next, we use the zero boundary and initial conditions to obtain the second loss as follows:
\[
    \mathcal{L}_{Operator}(\theta) = \frac{1}{NP}\sum_{i=1}^N\sum_{k=1}^P| G_{\theta}(g^i)(y_{g, k}^i) -G(g^i)(y_{g, k}^i)|^2
\]
where $y_{g, k}^i = (\tau_{g, k}^i, x_{g,k}^i)$ denotes points from the initial and boundary conditions. Hence, the total loss becomes
\begin{equation}
\label{totalloss}
    \mathcal{L}(\theta)=\mathcal{L}_{Physics} (\theta) + \mathcal{L}_{Operator}(\theta).
\end{equation}
It follows that by minimizing the loss function (\ref{totalloss}) the network can effectively predict the solution of the HJB equation. Fig. \ref{fig_piDeepONet} shows the schematics of physics-informed DeepONet connected in a feedforward manner.

\section{Results and Discussion}
\begin{figure}
   \begin{center}
           \includegraphics[width=\textwidth]{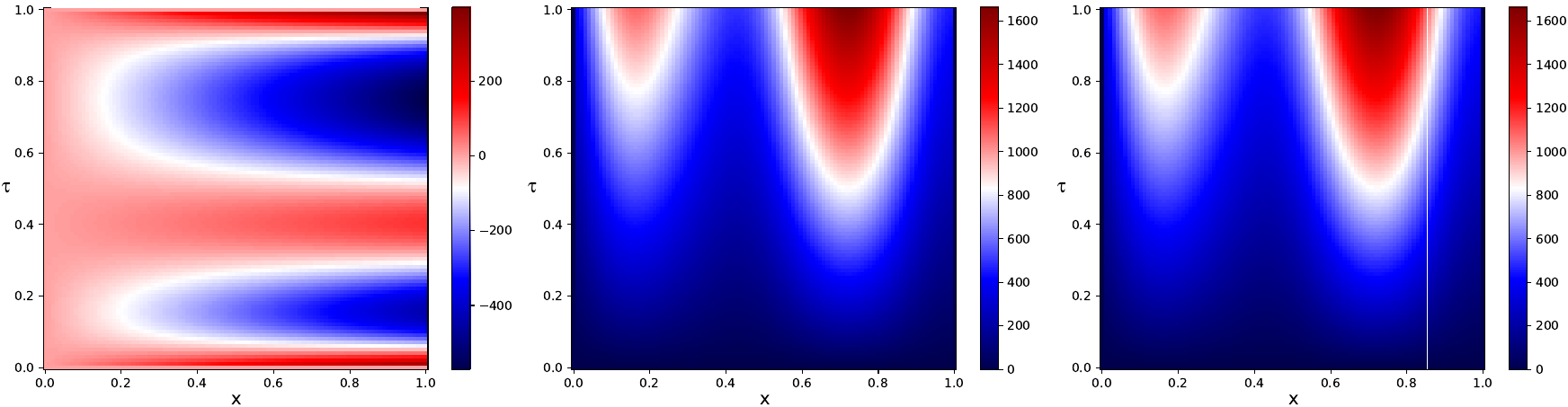}

{\small Input function $g$ \hglue2truecm NN solution \hglue1.5truecm FDM numerical solution}

   \end{center}

\caption{\small Comparision of a physics-informed DeepONet solution (NN solution) and the numerical solution obtained by the finite difference method (FDM numerical solution). The right-hand side represents the input function $g$.}

  \label{figNN_arch}
\end{figure}

The PI-DeepONet exhibits infinitesimal optimization and generalization errors, as it is easy to train and generalizes well to unseen data. In our approach, we did not use any input-output data, rather we only used the zero boundary and initial conditions. We approximate the PDE solution operator using branch and trunk nets. As a test example, we consider the diffusion function $\alpha(\varphi) = \varphi^2$. First, the input function of the branch net is discretized in a finite-dimensional space using a finite number of points called sensors. Then, the discretized input function is evaluated at fixed sensors to obtain point-wise evaluations. The trunk net takes the spatial and temporal coordinates and evaluates the solution operator to obtain the loss function. 
To generate our training data, we randomly sample $N = 500$ source term functions as input functions of the trunk net from a zero mean Gaussian process with an exponential quadratic kernel having a 0.2-length scale. The kernel function defines the covariance between two points in the process as a function of the distance between them. The parameter $l>0$ determines how quickly the covariance between two points decays as the distance between them increases. In this study, we set $l=0.2$. A smaller length scale results in a higher correlation between nearby points, whereas a larger length scale results in a lower correlation between nearby points.
Then, the selected input functions are evaluated at $m=100$ points as input sensors. The $m$ outputs of the source term functions are sent to the branch network. Next, we select the $P=100$ output sensors from the initial and boundary conditions, which are sent to the trunk nets. Our operator is then approximated by computing the dot product between the branch and trunk networks, and the corresponding operator loss is computed. After that, we select $Q=100$ collocation points inside the domain, and the error related to the underlying physics is computed. Finally, the total loss is evaluated by combining the two losses, which are minimized using the adaptive moment estimation (ADAM) optimizer with a learning rate of $10^{-3}$. Similarly, the test set is generated using the same approach. In Fig.~\ref{figNN_arch}, we compare a solution obtained by a physics-informed DeepONet method using the Relu activation function for 10000 iterations with a numerical solution constructed by means of the finite difference numerical method. 

\section{Conclusions}
In this study, we employed a physics-informed DeepONet to approximate the solution operator of a parametric parabolic equation arising from portfolio selection problems. The input function of the branch net was discretized in a finite-dimensional space using a fixed number of sensors. The discretized input functions were evaluated at fixed sensors to obtain point-wise evaluations. The operator was approximated by computing the dot product between the branch and trunk networks, and the corresponding operator loss is computed. We applied the physics-informed DeepONet to solve the model nonlinear parabolic equation obtained from the Hamilton-Jacobi-Bellman equation for solving optimal stochastic dynamic optimization problem.

\medskip
\noindent{\bf Acknowledgments.} The research was supported by the APVV-20-0311 (C.U.) and VEGA 1/0611/21 (D.\v{S}.) projects.

%


\begin{thebibliography}{6}
%


\bibitem{sirignano2018dgm}
Justin, S., Konstantinos, S.:
\newblock{DGM: A deep learning algorithm for solving partial differential equations}. 
\newblock {\em Journal of Comp. Physics},  375, 1339--1364 (2018).

\bibitem{KilianovaSevcovicANZIAM}
Kilianov\'a, S., \v{S}ev\v{c}ovi\v{c}, D.:
\newblock{A Transformation Method for Solving the Hamilton-Jacobi-Bellman Equation for a Constrained Dynamic Stochastic Optimal Allocation Problem}. 
\newblock {\em ANZIAM Journal}, 55, 14--38, (2013).


\bibitem{KilianovaTrnovska}
Kilianov\'a, S.,  Trnovsk\'a, M.:
\newblock {Robust Portfolio Optimization via solution to the Hamilton-Jacobi-Bellman Equation}. 
\newblock {\em Int. Journal of Comp. Math.}, 93, 725--734 (2016).


\bibitem{lu2021learning}
Lu, Lu, Jin, P., Pang, G., Zhang, Z., Karniadakis, G. Em.:
\newblock{Learning nonlinear operators via DeepONet based on the universal approximation theorem of operators}.
\newblock {\em Nature Machine Intelligence},  3, 218--229  (2021).


\bibitem{raissi2019physics}
Raissi, M., Perdikaris, P., Karniadakis, G. E.: 
\newblock{Physics-informed neural networks: A deep learning framework for solving forward and inverse problems involving nonlinear partial differential equations}. 
\newblock {\em Journal of Comp. Physics},  378, 686--707 (2019).

\bibitem{udeani2021application}
{\v{S}}ev{\v{c}}ovi{\v{c}}, D., Udeani, C.I.:
\newblock{Application of maximal monotone operator method for solving Hamilton--Jacobi--Bellman equation arising from optimal portfolio selection problem}.
\newblock {\em Jpn. J. Ind. Appl. Math.}, 5, pp 1--21 (2021).


\bibitem{chen1995universal}
Tianping, C., Hong, C.:
\newblock{Universal approximation to nonlinear operators by neural networks with arbitrary activation functions and its application to dynamical systems}.
\newblock {\em IEEE Transactions on Neural Networks}, 6, 911--917 (1995).


\bibitem{zhu2019physics}
Zhu, Y., Zabaras, N., Koutsourelakis, P., Perdikaris, P.: 
\newblock {Physics-constrained deep learning for high-dimensional surrogate modeling and uncertainty quantification without labeled data}. 
\newblock {\em Journal of Comp. Physics}, 394, 56--81 (2019).



\end{thebibliography}
\end{document}